\documentclass[11pt,a4]{article}
\usepackage{amsmath}
\usepackage{amssymb}
\usepackage{theorem}
\usepackage{color}
\newcommand{\R}{{\cal R}}
\newcommand{\B}{{\cal B}}

\newtheorem{theorem}{Theorem}

{\theorembodyfont{\rmfamily}
{\theorembodyfont{\rmfamily}

\newcommand{\scdot}{{\cdot\,}}
\newcommand{\svert}{{\vert\,}}

\begin{document}
{\small \hfill Probability Symposium, RIMS Workshop\\
\hfill Dec.~19-22, 2016 }
\begin{center}
Recent progress on conditional randomness\\
Hayato Takahashi\footnote{
hayato.takahashi@ieee.org\newline
This work was done when the author was with Gifu University and supported by KAKENHI (24540153).
}
\end{center}
\date{}

\thispagestyle{empty}
The set of Hippocratic random sequences w.r.t.~\(P\) is defined as the compliment of the effective null sets w.r.t.~\(P\) and denote it by \(\R^P\).
In particular if \(P\) is computable it is called Martin-L\"of random sequences. 

Lambalgen's theorem (1987) \cite{lambalgen87} says that a pair of sequences \((x^\infty,y^\infty)\in\Omega^2\) is Martin-L\"of (ML) random w.r.t.~the product  of uniform measures iff \(x^\infty\) is ML-random and \(y^\infty\) is ML-random relative to \(x^\infty\),
where \(\Omega\) is the set of infinite binary sequences. 
In \cite{{vovkandvyugin93},{takahashiISIT2006},{takahashiIandC},{takahashiIandC2}}, generalized Lambalgen's theorem is studied for computable correlated probabilities.

%
%


Let \(S\) be the set of finite binary strings and \(\Delta(s):=\{sx^\infty \vert x^\infty\in\Omega\}\) for \(s\in S\), where \(sx^\infty\) is the concatenation of \(s\) and \(x^\infty\).
Let \(X=Y=\Omega\) and  \(P\) be a computable probability on \(X\times Y\).  \(P_X\) and \( P_Y\) are marginal distribution on \(X\) and \(Y\), respectively. 
In the following we write \(P(x,y):=P(\Delta(x)\times\Delta(y))\) and \(P(x\vert y):=P(\Delta(x)\vert\Delta(y))\) for \(x,y\in S\).

Let \(\R^P\) be the set of ML-random points and \(\R^P_{y^\infty}:=\{ x^\infty\mid (x^\infty,y^\infty)\in\R^P\}\).
In \cite{{takahashiISIT2006},{takahashiIandC}}, it is shown that conditional probabilities exist for all random parameters, i.e., 
\[
\forall x\in S,\ y^\infty\in\R^{P_Y}\ \ P(x\vert y^\infty):=\lim_{y\to y^\infty}P(x\vert y)\text{ (the right-hand-side exist)}
\]
and \(P(\cdot\vert y^\infty)\) is a probability on \((\Omega, \B)\) for each \(y^\infty\in\R^{P_Y}\).

Let \(\R^{P(\cdot\vert y^\infty),y^\infty}\) be the set of Hippocratic random sequences w.r.t.~\(P(\cdot\vert y^\infty)\) with oracle \(y^\infty\).

\begin{theorem}[\cite{{takahashiISIT2006},{takahashiIandC},{takahashiIandC2}}]\label{col-ifpart}
Let \(P\) be a computable probability on \(X\times Y\). Then 
\begin{equation}
\R^P_{y^\infty}\supseteq\R^{P(\cdot\vert y^\infty),y^\infty} \text{ for all }y^\infty\in\R^{P_Y}.
\end{equation}
Fix \(y^\infty\in\R^{P_Y}\) and suppose that \(P(\cdot\vert y^\infty)\) is computable with oracle \(y^\infty\).
Then 
\begin{equation}\label{generalLambal}
\R^P_{y^\infty}=\R^{P(\cdot\vert y^\infty),y^\infty}.
\end{equation}
\end{theorem}

\newpage
\thispagestyle{empty}
It is known that there is a non-computable conditional probabilities \cite{roy2011} and in \cite{Bauwens2015} Bauwens showed an example that violates the equality in (\ref{generalLambal})
when the conditional probability is not computable with oracle \(y^\infty\). In \cite{takahashi2014}, an example that  for all \(y^\infty\), the  conditional probabilities are not computable with oracle \(y^\infty\) and (\ref{generalLambal}) holds.
A survey on the randomness for conditional probabilities  is shown in \cite{BST2016}.

%
%
%

Next we study mutually singular conditional probabilities.
In \cite{kjoshanssen},  Hanssen showed that for Bernoulli model \(P(\cdot\vert \theta)\),
\begin{equation}\label{hanssen}
\R^{P(\cdot\vert\theta)}=\R^{P(\cdot\vert\theta),\theta}\text{ for all }\theta.
\end{equation}
We generalize  Hanssen's theorem (\ref{hanssen}) for mutually singular conditional  probabilities. 
In \cite{{takahashiISIT2006},{takahashiIandC2}}, equivalent conditions for mutually singular conditional probabilities are shown. 
\begin{theorem}[\cite{{takahashiISIT2006},{takahashiIandC2}}]\label{th-consis-pos}
Let \(P\) be a computable probability on \(X\times Y\), where \(X=Y=\Omega\).
The following six statements are equivalent:\\
(1) \( P(\scdot\svert  y)\perp  P(\scdot\svert  z)\) if \(\Delta( y)\cap\Delta( z)=\emptyset\) for \( y, z\in S\).\\
(2) \(\R^{ P(\scdot\svert  y)}\cap\R^{ P(\scdot\svert  z)}=\emptyset\) if \(\Delta( y)\cap\Delta( z)=\emptyset\) for \( y, z\in S\).\\
(3)  \( P_{Y|X}(\scdot\svert x)\) converges weakly to \(I_{ y^\infty}\) as \(x\to x^\infty\) for \((x^\infty, y^\infty)\in\R^{ P}\), where \(I_{ y^\infty}\) is the probability that has probability of 1 at \( y^\infty\).\\
(4) \(\R^{ P}_{ y^\infty}\cap\R^{ P}_ {z^\infty}=\emptyset\) if \( y^\infty\ne z^\infty\).\\
(5) There exists \(f:X\to Y\) such that  \(f(x^\infty)= y^\infty\) for \((x^\infty, y^\infty)\in\R^{ P}\).\\
(6) There exists \(f:X\to Y\) and \( Y'\subset Y\) such that \(P_Y( Y')=1\) and \(f= y^\infty,\ P(\scdot;  y^\infty)-a.s.\) for \( y^\infty\in Y'\).
\end{theorem}

Generalized form of Hanssen's theorem (\ref{hanssen}) is as follows.
\begin{theorem}
Let \(P\) be a computable probability on \(X\times Y\), where \(X=Y=\Omega\).
Under one of the condition of Theorem~\ref{th-consis-pos}, we have 
\[\R^P_{y^\infty}\supseteq\R^{P(\cdot\vert y^\infty)} \text{ for all }y^\infty\in\R^{P_Y}.\]
Fix \(y^\infty\in\R^{P_Y}\) and suppose that \(P(\cdot\vert y^\infty)\) is computable with oracle \(y^\infty\).
Then
\[\R^P_{y^\infty}=\R^{P(\cdot\vert y^\infty)} =\R^{P(\cdot\vert y^\infty),y^\infty}.\]
\end{theorem}

\newpage
\thispagestyle{empty}
{\small
\bibliographystyle{plain}

}

\end{document}